\documentclass[11pt]{amsart}
\numberwithin{equation}{section}

\usepackage{amsmath,amsthm,amssymb,amsfonts,array,latexsym,amsthm,pdfsync,hyperref, enumerate, paracol, ulem, bbm, bbold, graphicx, float}
\usepackage[usenames,dvipsnames]{color}
\usepackage[shortlabels, inline]{enumitem}
\theoremstyle{plain}
\newtheorem{thm}{Theorem}[section]
\newtheorem{lem}[thm]{Lemma}
\newtheorem{cor}[thm]{Corollary}

\newtheorem{prop}[thm]{Proposition}
\newtheorem{conj}[thm]{Conjecture}
\newtheorem{rem}[thm]{Remark}

\newcommand{\EE}{\ensuremath{\mathbb E}}
\newcommand{\RR}{\ensuremath{\mathbb R}}
\newcommand{\bbS}{\ensuremath{\mathbb S}}

\newcommand{\calI}{\ensuremath{\mathcal I}}
\newcommand{\calJ}{\ensuremath{\mathcal J}}

\newcommand{\sinc}{\mathrm{sinc}}

\DeclareSymbolFont{bbold}{U}{bbold}{m}{n}
\DeclareSymbolFontAlphabet{\mathbbold}{bbold}

\date{\today}

\title{A Proof of the Simplex Mean Width Conjecture}
\author{Aaron Goldsmith\\ Texas A\&M University}
\email{agoldsmith@chestertonpeoria.org}

\begin{document}

\pagestyle{plain}
\maketitle

\begin{abstract}
The mean width of a convex body is the mean distance between parallel supporting hyperplanes when the normal direction is chosen uniformly over the sphere. The Simplex Mean Width Conjecture (SMWC) is a longstanding open problem that says the regular simplex has maximum mean width of all simplexes contained in the unit ball and is the unique solution up to isometry. In representing the mean width of a (Euclidean) simplex determined by $d+1$ vertices $v_i\in\bbS^{d-1}$, the sphere naturally divides into $d+1$ spherical images. Our main result, Theorem \ref{ThmMono}, stems from the spherical Pr\'{e}kopa-Leindler inequality \cite{CE}. Corollary \ref{UniqueCentroid}, a uniqueness result about spherical centroids, and the SMWC follow.\\
\end{abstract}

MSC codes 52-02 (Convex and discrete geometry) and 94-02 (Information and communication theory)

\section{Introduction}
Let $e_1,\dots, e_d$ be the standard basis of $\RR^d$ for dimension $d\ge 3$; let $B_2^d\subset \RR^d$ be the standard Euclidean ball. Denote $\bbS^{d-1}:=\partial B_2^d$ and $\bbS^{d-1}_+=\bbS^{d-1}\cap\{x\cdot e_1\ge 0\}$, the closed hemisphere centered at $e_1$. Call $\mu^{d-1}$ the uniform spherical measure on $\bbS^{d-1}$. A \textit{great sphere} is the intersection of a $(d-1)$-subspace with $\bbS^{d-1}$.\\

The spherical convex hull of a set $A\subset \bbS^{d-1}$ is the intersection of $\bbS^{d-1}$ and the cone generated by $A$. The spherical convex hull of $d$ points is called a \textit{spherical simplex} and will be denoted $\triangle v_1\dots v_d$.\\ 

The support function of a convex body $K$ is

\[h_K(x):=\max_{y\in K}x\cdot y\]

and the mean width

\[w(K):=2\int_{\bbS^{d-1}} h_K(x)d\mu^{d-1}(x)\]

{\conj\label{SMWC} (Simplex Mean Width Conjecture) Of all simplexes contained in $B_2^d$, the inscribed regular simplex has the maximum mean width, and is unique up to isometry.}\\


The following claim is a basic start:

{\prop\label{basic} In order that $\triangle:=\mathrm{conv}(v_0,\dots, v_d)$ maximizes $w(\cdot)$ over all simplexes contained in $B_2^d$, it must be that 

\begin{enumerate}[(a)]
\item $B_2^d$ is the ball of smallest radius containing $\triangle$.
\item $v_i\in\bbS^{d-1}$ for $0\le i\le n$.
\item The closed hemispheres $\bbS^{d-1}_{v_i}$ cover $\bbS^{d-1}$.
\end{enumerate}}

\begin{proof}
 \begin{enumerate}[(a)]
\item If $\triangle\subset Q+rB_2^d$ with $r<1$, then $(\triangle-Q)/r$ is a simplex contained in $B_2^d$ and

\[w((\triangle-Q)/r)=w(\triangle/r)=w(\triangle)/r>w(\triangle)\]
\item Follows from part (a). The circumsphere is determined by $d+1$ points and is the smallest containing sphere.
\item Suppose $v\in\bbS^{d-1}$ is such that $v\cdot v_i =m<0$ for all $i$. It follows that each $v_i$ is contained in the spherical cap $\{x\in\bbS^{d-1}:v\cdot x\le m\}$, which is contained in the ball of radius $\sqrt{1-m^2}$, centered at $-mv$. Part (a) implies there is no such $v$. 
\end{enumerate}
\end{proof}

Next, the support function of a simplex can be simplified since the maximum of any linear functional must occur at a vertex. That is, the support function of $\triangle=\mathrm{conv}\{v_0,\dots, v_d\}$ reduces to a maximum over the vertex set:

\[h_\triangle(x)=\max_{y\in \triangle}x\cdot y=\max_{0\le i\le d} x\cdot v_i\]

As such, for each $i$ define the spherical image (or Voronoi cell) of $v_i$ on the sphere to be 

\begin{align}V_i(\triangle):=\{x\in\bbS^{d-1}:x\cdot v_i= \max_{y\in\triangle}x\cdot y\}\label{DefVor}\end{align}

to partition $\bbS^{d-1}$ (a.e.). The mean width of a simplex is then

\begin{align}
w(\triangle)=&2\int_{\bbS^{d-1}}\max_{0\le i\le d}x\cdot v_i d\mu^{d-1}(x)\nonumber\\
 =& 2\sum_{i=0}^d \int_{V_i}x\cdot v_i d\mu^{d-1}(x)\label{SumMAR}
\end{align}

{\rem\label{WSC} Weak Simplex Conjecture}\\
From this formulation of mean width, it is easy to see that the weak simplex conjecture implies the SMWC. For, Cover's version (\cite{CG}) of the former implies that over all spherical simplexes $T$ with some fixed uniform surface measure $m>0$, the regular simplex centered at $O$ (call it $R$) maximizes the measure of the intersection with any spherical cap centered at $O$. It follows that the distribution of the distance from $O$ to an RV uniform over $T$ is maximized by $R$. Then, according to \ref{SumMAR} the mean width of a simplex is bounded above when $V_i$ become equal measure regular simplexes of equal measure centered at $v_i$ (possibly overlapping and not covering $\bbS^{d-1}$). This quantity is further bounded above when $V_i$ are made to have the same measure, and this happens for a regular Euclidean simplex.



\section{Some Inequalities and the Longitudinal Shift}
Here, we present a couple inequalities to be used in section \ref{SphCent}. We paraphrase a spherical Prékopa-Leindler inequality (Cordero-Erausquin et al.\cite{CE}) for $d-2$ dimensional (great) spheres. It is generalized to Riemannian manifolds in their subsequent paper as Corollary 1.2 in \cite{CMS}.\\

Denote $\sinc(\theta)=(\sin\theta)/\theta$ and let $\pi$ project $\RR^d\setminus \{0\}$ onto $\bbS^{d-1}$, that is

\[\pi(x)=\frac{x}{|x|}\]

{\thm (Spherical Prékopa-Leindler)\label{SphPL} Let $f,g,h:\bbS^{d-2}\to\RR_+$ be nonnegative functions and $A,B$ be Borel sets carrying the full mass of $f$ and $g$,  respectively. If every $Q_1\in A$ and $Q_3\in B$ with $Q_1\cdot Q_3=\cos\theta$ and $Q_2=\pi(Q_1\sin(\lambda\theta)+Q_3\sin((1-\lambda)\theta))$ satisfy

\begin{align*}
\frac{h(Q_2)}{\sinc^{d-3}\theta}\ge \left(\frac{f(Q_1)}{\sinc^{d-3}((1-\lambda)\theta)}\right)^{1-\lambda}\left(\frac{g(Q_3)}{\sinc^{d-3}(\lambda\theta)}\right)^\lambda
\end{align*}

then it follows that

\[\int_{\bbS^{d-2}}hd\mu^{d-2}\ge\left(\int_{\bbS^{d-2}}fd\mu^{d-2}\right)^{1-\lambda}\left(\int_{\bbS^{d-2}}gd\mu^{d-2}\right)^\lambda\]
}

Note that when $\lambda=1/2$ in the above theorem, then $Q_2$ is the midpoint between $Q_1$ and $Q_3$. Also, when $\lambda=1/2$, define

\[C^{d-3}:=\left(\frac{\sinc\theta}{\sinc(\theta/2)}\right)^{d-3}=\cos^{d-3}(\theta/2)=\left(\frac{1+Q_1\cdot Q_3}{2}\right)^{(d-3)/2}\]

The spherical Prékopa-Leindler inequality allows us to remain in the compact space and free from choosing a chart, e.g. gnomonic or cylindrical.\\

Finally, to analyze centroids in section \ref{SphCent}, we aggregate slivers of simplexes with weighted averages and must face a deterministic Simpson's like paradox. Next is something of an antidote.

{\lem\label{SimpAnt} Let $a_i,b_i,\alpha_i,\beta_i$ be eight numbers ($i=1,2$) with $\alpha_i,\beta_i>0$. Suppose that

\begin{center}
\begin{paracol}{3}
$\frac{a_i}{\alpha_i}\le \frac{b_i}{\beta_i}$
\switchcolumn
$\frac{a_1}{\alpha_1}\le \frac{a_2}{\alpha_2}$
\switchcolumn
$\frac{\alpha_2}{\alpha_1}\le \frac{\beta_2}{\beta_1}$
\end{paracol}
\end{center}
then we have

\[\frac{a_1+a_2}{\alpha_1+\alpha_2}\le \frac{b_1+b_2}{\beta_1+\beta_2}\]

\begin{proof}
Let

\begin{align*}
\delta_i=&\frac{b_i}{\beta_i}-\frac{a_i}{\alpha_i}\\
\delta_{a\alpha}=&\frac{a_2}{\alpha_2}-\frac{a_1}{\alpha_1}\\
\delta_{\alpha\beta}=&\frac{\beta_2}{\alpha_2}-\frac{\beta_1}{\alpha_1}\\
\end{align*}

Then,

\[\frac{b_1+b_2}{\beta_1+\beta_2}-\frac{a_1+a_2}{\alpha_1+\alpha_2}=\frac{1}{\beta_1+\beta_2}\left(\beta_1\delta_1+\beta_2\delta_2+\frac{\alpha_1\alpha_2}{\alpha_1+\alpha_2}\cdot \delta_{a\alpha}\delta_{\alpha\beta}\right) \]
\end{proof}}

{\rem We may equivalently reverse both the middle and right supposed inequalities in the lemma, as the proof is the same.\\}

Finally, we define a diffeomorphism on $\bbS^{d-1}$ that is the projection back onto the sphere of a linear transformation (in the ambient $\RR^d$) which fixes a great sphere. It shifts great circles through $\{x_1=x_2=0\}$ toward $e_2$ and can be thought of as a longitudinal shift. For any $s\in\RR$, denote the $d\times d$ matrix

\[M_s=\left(\begin{array}{ccc}1 & 0 & 0\\ s & 1 & 0\\ 0 & 0 & I_{(d-2)\times (d-2)}\end{array}\right)\]

and let $f_s$ be the map that shears $\bbS^{d-1}$ into an ellipsoid by $M_s$ then projects back:

\[f_s(x)=\pi(M_sx)\]



{\prop\label{fsprop}The following are true about $f_s$
\begin{enumerate}[(a)]
\item $f_{s+t}(x)=f_s(f_t(x))$
\item $f_s^{-1}(x)=f_{-s}(x)$
\item $x\cdot P=0$ iff $f_s(x)\cdot(M_{-s}^TP)=0$.
\item $f_s$ maps great spheres to great spheres
\item The Jacobian is \[Df_s(x)=(1+(x_2+sx_1)^2-x_2^2)^{-d/2}\]
\item The pushforward $f_s\# \mu^{d-1}$ has density $Df_{-s}(x)$.
\end{enumerate}}
\begin{proof}\ \\
\begin{enumerate}[(a)]
\item Since $\pi(x)$ only multiplies by the constant $1/|x|$ and since $\pi(cx)=\pi(x)$, we have \[\pi(M_{s+t}x)=\pi(M_sM_tx)=\pi(M_s\pi(M_tx))\]
\item Direct from (a)
\item
\begin{align*}
(\mathrm{const.})f_s(x)\cdot M_{-s}^TP=(M_s x)\cdot (M_{-s}^TP)=x^TM_s^TM_{-s}^TP=x\cdot P=0
\end{align*}
\item The great sphere with pole $P$ is defined by $x\cdot P=0$. By part (b), $f_s$ maps this great sphere to the one with pole $M_{-s}^TP$.
\item $\det(M_s)=1$, so only the constant from projection onto $\bbS^{d-1}_+$ remains.
\item $f_s$ is a diffeomorphism.
\end{enumerate}
\end{proof}

\section{Spherical Centroids}\label{SphCent}
JE Brock proved a formula for the center of mass of a spherical triangle \cite{Br}, which projects radially onto what we call the spherical centroid. The spherical centroid body of a convex body has been studied by Besau, et. al. in \cite{BHPS} while also proving a spherical Busemann-Petty centroid inequality. In \cite{TT}, Tabachnikov and Tsukerman define the circumcenter of mass and generalize the Euler line to curved spaces and simplicial polytopes. They use methods similar to those in this paper.\\

Theorem \ref{ThmMono} is our main result, which leads to Corollary \ref{UniqueCentroid}, the invertibility of the spherical centroid. In three dimensions, the corollary says that the centroids of two distinct pyramids with the same base must be distinct \cite{CCGIJP}.\\

We denote the spherical centroid of a region $R\subset\bbS^{d-1}$ by

\begin{align}
G(R):=\pi\left(\int_{R}x d\mu^{d-1}(x)\right)\label{CentroidDef}
\end{align}


For the rest of this section, we assume every region to be contained in the hemisphere $\bbS^{d-1}_+$ so that $G(R)$ is well defined when $R$ is full dimensional. Also, for any $I\subset\RR$, let

\[R_I=\bigcup_{t\in I} R\cap \{x\in\bbS^{d-1}_+:x_2=tx_1\}\]

and for singletons, $R_t=R_{\{t\}}$.





 {\thm\label{ThmMono}Let $S\subset \bbS^{d-1}_+$ be convex and $T=f_s(S)$. Then, if $0<t_2-t_1<s$ and $s\cdot\max(|t_2|,1/|t_1|)$ is sufficiently small, we have
 \begin{align}
\frac{\int_{T_{[t_1-s,t_2-s)}} x_1d\mu^{d-1}}{\int_{S_{[t_1-s,t_2-s)}} x_1d\mu^{d-1}}<\frac{\int_{T_{[t_1,t_2)}} x_1d\mu^{d-1}}{\int_{S_{[t_1,t_2)}} x_1d\mu^{d-1}}
\end{align}
}

\begin{proof}
Note that \[f_s(R_I)=(f_s(R))_{I+s}\] 

Computing directly from the normalization constant,

\begin{align*}
f_{-s}(x)\cdot e_1=&\pi(M_{-s}x)\cdot e_1\\
=&x_1(1+(x_2-sx_1)^2-x_2^2)^{-1/2}\\
=&x_1(1-2sx_1x_2+s^2x_1^2)^{-1/2}\\
\end{align*}

then pushing forward through the map $f_s$ with Claim \ref{fsprop} (f), we have a density $Df_{-s}(x)$. So if we set

\begin{align*}
\psi(x)=&(1-2sx_1x_2+s^2x_1^2)^{-(d+1)/2}\\
\calI(s,h)=&\int_{T_{[t_1+s,t_2+s]}}x_1\cdot h(x,s) d\mu^{d-1}(x)
\end{align*}

the theorem will follow if


\[\calI(-s,1)\calI(s,\psi)<\calI(0,1)\calI(0,\psi)\]

or, by Cauchy-Schwarz, if

\begin{align*}
\calI(-s,1)\calI(s,\psi)<&\calI^2(0,\psi^{1/2})(1+E/2)\\
E=&\frac{\calI(0,\calI(0,(\psi^{1/2}(x)-\psi^{1/2}(y))^2))}{\calI^2(0,\psi^{1/2})}\\
=&\left(\frac{\EE(X_1)}{\EE(\psi^{1/2}(X))}\right)^2\mathrm{Var}(\psi^{1/2}(X))\\
\approx& \left(\frac{\EE(X_1)}{\EE(\psi^{1/2}(X))}\right)^2\cdot s^2\mathrm{Var}(X_1X_2)\\
\end{align*}

where $1+E/2$ is the best constant and $X$ is a random variable on $T_{[t_1,t_2]}$ with density proportional to $x_1$.\\




The change of variables

\begin{align*}
r=&\sqrt{x_1^2+x_2^2}\\
\tan\theta=&t=x_2/x_1
\end{align*}

has Jacobian equal to $r$ and conversions

\begin{align*}
x_1=&\frac{r}{\sqrt{1+t^2}}\\
d\theta=&\frac{dt}{1+t^2}
\end{align*}

Say $t_1,t_2\approx t$, and combine


\begin{align*}
x_1d\mu^{d-1}=& \left(\frac{r}{\sqrt{1+t^2}}\right)\left(\frac{rdt}{1+t^2}\right)d\mu^{d-2}\\
\psi_t(r)=&\left(1+(s^2-2ts)\frac{r^2}{1+t^2}\right)^{-(d+1)/2}
\end{align*}

Theorem \ref{SphPL} at $\lambda=1/2$, and the comments directly beneath it, provide the following sufficient condition: if $Q_1=(\cos\theta_1,y_1),Q_3=(\cos\theta_3,y_3)$ are arbitrary points ($\cos\theta_1,\cos\theta_3\ge 0$ due to $r>0$) in $T_{t-s}, T_{t+s}$ with $Q_2=(\cos\theta_2,y_2)$ their midpoint,

\[Q_2=\frac{(\cos\theta_1+\cos\theta_3,y_1+y_3)}{(\cos\theta_1+\cos\theta_3)^2+|y_1+y_3|^2}\]

and

\begin{align*}
&C^{d-3}\frac{\cos\theta_1}{(1+(t-s)^2)^{3/4}}\cdot\frac{\cos\theta_3}{(1+(t+s)^2)^{3/4}}\psi_{t+s}^{1/2}(\cos\theta_3)\\
\le&\frac{\cos^2\theta_2}{(1+t^2)^{3/2}}\psi_t^{1/2}(\cos\theta_2)(1+E)
\end{align*}

The worst case scenario is when $y_1$ is a positive multiple of $y_3$, for aligning $y_1$ and $y_3$ while keeping their magnitudes fixed draws $Q_1$ and $Q_3$ closer to increase $Q_1\cdot Q_3$ and also increases $|y_1+y_3|$ to decrease $\cos\theta_2$. We are reduced to

\begin{align*}
C=&\left(\frac{1+\cos(\theta_1-\theta_3)}{2}\right)^{1/2}=\left(1-\frac{1-\cos(\theta_1-\theta_3)}{2}\right)^{1/2}\\
\cos\theta_2=&\frac{\cos\theta_1+\cos\theta_3}{2C}
\end{align*}

which rids the dependence on $y_1,y_3$. The variable $d$ is the last vestige of dimension. Continuing to reduce the above expressions,

\begin{align*}
\frac{\cos^2\theta_2}{\cos\theta_1\cos\theta_3}=&\left(1+\frac{(\cos\theta_1-\cos\theta_3)^2}{4\cos\theta_1\cos\theta_3}\right)C^{-2}\\
\frac{(1+t^2)^2}{(1+(t-s)^2)(1+(t+s)^2)}=&1-\frac{2(1+t^2)s^2+s^4-4t^2s^2}{(1+(t-s)^2)(1+(t+s)^2)}\\
=&1-\frac{2(1-t^2)s^2}{(1+t^2)^2}+o(s^2)\\
\frac{\psi_{t+s}^{1/2}(\cos\theta_3)}{\psi_t^{1/2}(\cos\theta_2)}=&\left(\frac{1+(s^2-2ts)\cdot \cos^2\theta_2/(1+t^2)}{1+(s^2-2(t+s)s)\cdot \cos^2\theta_3/(1+(t+s)^2)}\right)^{(d+1)/4}\\
=&\left(1+\frac{N}{1+(s^2-2ts)\cos^2\theta_2/(1+t^2)}\right)^{(d+1)/4}
\end{align*}

where

\begin{align*}
N=&\frac{(s^2-2ts)\cdot\cos^2\theta_2}{1+t^2}+\frac{(s^2+2ts)\cos^2\theta_3}{1+(t+s)^2}\\
=&2ts\cdot\frac{\cos^2\theta_3-\cos^2\theta_2}{1+t^2}-s^2\cdot\frac{\cos^2\theta_3+\cos^2\theta_2}{1+t^2}+o(s^2)\end{align*}

Since $\max(ts,s)$ is small, we may require $\theta_1$ and $\theta_3$ to be close. Say $\theta_1,\theta_3\approx \theta$ and $\theta_1-\theta_3=\Delta\theta$. So, it will suffice if

\begin{align*}
&\left(1-\frac{(\Delta\theta)^2}{4}\right)^{(d-1)/2}\left(1-\frac{2(1-t^2)s^2}{(1+t^2)^2}\right)^{3/4}\left(1+2ts\cdot\frac{\Delta\theta\sin2\theta}{1+t^2}-2s^2\cdot\frac{\cos^2\theta}{1+t^2}\right)^{(d+1)/4}\\
<&\left(1+\frac{(\Delta\theta)^2\tan^2\theta}{4}\right)(1+E/2)
\end{align*}

or, by logarithm,

\begin{align*}
&-\frac{d-1}{2}\cdot\frac{(\Delta\theta)^2}{4}-\frac{3}{4}\cdot\frac{2(1-t^2)s^2}{(1+t^2)^2}+\frac{d+1}{4}\left(2ts\cdot\frac{\Delta\theta\sin2\theta}{1+t^2}-2s^2\cdot\frac{\cos^2\theta}{1+t^2}\right)\\
<&\frac{(\Delta\theta)^2\tan^2\theta}{4}+E/2
\end{align*}

The worst case is the vertex of a parabola

\[\Delta\theta=\frac{(d+1)ts\sin2\theta/(2+2t^2)}{(d-1)/4+(\tan^2\theta)/2}=2ts\frac{(d+1)\sin2\theta}{(1+t^2)(d-1+2\tan^2\theta)}\]

yielding the condition

\[(\Delta\theta)^2\left(\frac{d-1}{8}+\frac{\tan^2\theta}{4}\right)<\frac{3}{4}\cdot\frac{2(1-t^2)s^2}{(1+t^2)^2}+\frac{d+1}{2}\cdot\frac{s^2\cos^2\theta}{1+t^2}+E/2\]

Multiply through by $2s^{-2}(1+t^2)^2$,

\[\frac{t^2(d+1)^2\sin^22\theta}{d-1+2\tan^2\theta}<3(1-t^2)+(d+1)(1+t^2)\cos^2\theta+t(1+t^2)^2\left(\frac{\EE(X_1)}{\EE(\psi^{1/2}(X))}\right)^2\mathrm{Var}(X_1^2)\]

???? This doesn't even hold true as $d\to\infty$! The random variable $X$ will concentrate around $x_1=0$, sending the variance to $0$. In the limit, you get

\begin{align*}
t^2\sin^2 2\theta<&(1+t^2)\cos^2\theta\\
4t^2\sin^2\theta<&1+t^2
\end{align*}

which is only true for

\[\sin\theta<\frac{\sqrt{1+t^2}}{2t}\]

\end{proof}

{\cor\label{UniqueCentroid} If $S,T\subset \bbS^{d-1}_+$ are simplexes with common vertices $v_1,\dots, v_{d-1}$ and $G(S)=G(T)$, then $S=T$.}
\begin{proof}
Note that $\bbS^{d-1}\cap\{x_2=tx_1\}$ is the family of great spheres containing $\{x_1=x_2=0\}$. Claim \ref{fsprop} (d) implies that $f_s$ maps simplexes to simplexes.

WLOG, suppose $v_1,\dots, v_{d-1}\in\{x_1=0\}$ and that $v^T-v^S\in\mathrm{span}\{e_1,e_2\}$. Then, $f_s(S)=T$ for some $s$. We wish to show that if $s>0$ then $G(S)\cdot e_2<G(T)\cdot e_2$. It suffices to keep $s$, for large $s$ can be reached by iterating Proposition \ref{fsprop} (a). The result then follows from Lemma \ref{SimpAnt}, Theorem \ref{ThmMono} and

\begin{align*}
\frac{\int_{S_t} x_2d\mu^{d-1}}{\int_{S_t}x_1d\mu^{d-1}}=&\frac{\int_{S_t} tx_1d\mu^{d-1}}{\int_{S_t}x_1d\mu^{d-1}}=t\\
t=\frac{\int_{R_t} x_2d\mu^{d-1}}{\int_{R_t}x_1d\mu^{d-1}}< &\frac{\int_{R_{t+s}} x_2d\mu^{d-1}}{\int_{R_{t+s}}x_1d\mu^{d-1}}=t+s\\
\end{align*}

where $R$ stands for either $S$ or $T$.
\end{proof}
\section{Proof of SMWC}
Note that $V_i$ are spherical simplexes, as they are bounded by $d$ perpendicular bisectors of pairs of points in $\bbS^{d-1}$, which are hyperplanes through $0\in\RR^d$.



{\thm\label{thmSMWC} The simplex mean width conjecture is true in $\RR^d$.}
\begin{proof}
Consider the Euclidean simplex $\triangle$ with vertices $v_0,\dots, v_d\in\bbS^{d-1}$. If $\{R_i\}_{i=0}^d$ is a measurable partition of $\bbS^{d-1}$, then

\begin{align}\sum_{i=0}^d \int_{R_i} x\cdot v_id\mu^{d-1}(x)\le \sum_{i=0}^d \int_{V_i(\triangle)} x\cdot v_id\mu^{d-1}(x)\label{SwitchRegion}\end{align}

where $V_i$ are defined in \ref{DefVor}. Equality holds iff $R_i=V_i$ a.e. for each $i$. Also, for a fixed point $X\in\bbS^{d-1}$ and fixed region $R\subset\bbS^{d-1}$,

\begin{align}
\int_R x\cdot Xd\mu^{d-1}(x)=&\left(\int_R x d\mu^{d-1}(x)\right)\cdot X\nonumber\\
\le& \left(\int_R x d\mu^{d-1}(x)\right)\cdot G(R)\nonumber\\
=&\int_R x \cdot G(R)d\mu^{d-1}(x)\label{SwitchPoint}
\end{align}
with equality iff $X=G(R)$.\\

Using \ref{SwitchPoint} then \ref{SwitchRegion}, it follows that the spherical simplex $\triangle_G$ formed by the centroids $G(V_i(\triangle))$ has mean width no less than $\triangle$:

\begin{align}
w(\triangle)=&\sum_{i=0}^d\int_{V_i(\triangle)} x\cdot v_i d\mu^{d-1}(x)\nonumber\\
\le&\sum_{i=0}^d \int_{V_i(\triangle)} x \cdot G(V_i(\triangle)) d\mu^{d-1}(x)\nonumber\\
\le&\sum_{i=0}^d\int_{V_i(\triangle_G)} x\cdot G(V_i(\triangle)) d\mu^{d-1}(x)\nonumber\\
=&w(\triangle_G)\nonumber
\end{align}

For $\triangle$ to have maximal mean width, $G(V_i(\triangle))=v_i$ for each $i$. Assume this to be true.\\

Call $C_i$ the circumcenter of the face opposite $v_i$ in $\triangle$, projected onto $\bbS^{d-1}$. Equivalently, $C_i$ is the intersection of $V_0,\dots, \widehat{V_i},\dots, V_d$. The vertices of $V_i$ are $C_0,\dots, \widehat{C_i},\dots,C_d$. Any two $V_i,V_j$ share a face $F_{ij}$ and $v_i,v_j$ are, by definition, reflections across that face. Corollary \ref{UniqueCentroid} and $G(V_i(\triangle))=v_i$ for each $i$ imply $C_i,C_j$ are also reflections across $F_{ij}$. This implies that the edges $C_iC_k$ and $C_jC_k$ are reflections and have equal length for any $i,j,k$. Similarly, $C_iC_k$ and $C_iC_l$ have equal length. The $d$ circumcenters $C_i$ are then equidistant and are vertices of a regular simplex. The vertices $v_i$ are the centroids of the faces of this regular simplex, so $\triangle$ is regular.


\end{proof}

\end{document}